\documentclass[11pt]{amsart}
\input{pramod.tex}
\input{tableau.tex}

\newcommand{\cF}{\mathcal{F}}
\newcommand{\bB}{\mathrm{B}}
\newcommand{\dd}{\mathbf{d}}
\newcommand{\gred}{G_\mathrm{red}}

\newcommand{\nilrad}{\mathfrak{n}}
\newcommand{\sz}[1]{\|#1\|^2}
\newcommand{\bN}{\nilcone_{o, r}}
\newcommand{\g}{\mathfrak{g}}
\newcommand{\bo}{\mathfrak{b}}
\newcommand{\h}{\mathfrak{h}}
\newcommand{\dualnil}{\Langdual\nilcone}
\newcommand{\local}{\Weights^{loc}_{+}}

\title{Local Systems on Nilpotent Orbits and Weighted Dynkin Diagrams}

\author{Pramod N. Achar}
\address{University of Chicago\\
  Chicago, IL 60637}
\email{pramod@math.uchicago.edu}
\author{Eric N. Sommers}
\address{University of Massachusetts---Amherst\\
  Amherst, MA 01003}
\email{esommers@math.umass.edu}

\date{14 December 2001}

\begin{document}

\maketitle

\section{Introduction}

Let $G$ be a reductive algebraic group over the complex numbers, $B$ a
Borel subgroup of $G$, and $T$ a maximal torus of $B$. We denote by
$\Weights = \Weights(G)$ the weight lattice of $G$ with respect to
$T$, and by $\DomWeights = \DomWeights(G)$ the set of dominant weights
with respect to the positive roots defined by $B$. Let $\g$ be the Lie
algebra of $G$, and let $\nilcone$ denote the nilpotent cone in~$\Lie{g}$.

Now, let $e \in \nilcone$ be a nilpotent element, and let $\orb_e$ be
the orbit of $e$ in $\g$ under the adjoint action of $G$.  We write
$G^e$ for the centralizer of $e$ in $G$.  Let $\nilcone_{o}$ denote
the set of nilpotent orbits in $\g$, and $\bN$ the set of
$G$-conjugacy classes of pairs
\[
\{(e, \tau) \mid \text{$e \in \nilcone$ and $\tau$ an irreducible
rational representation of $G^e$} \}.
\]

Lusztig \cite{lusztig:cells-aff-4} conjectured the existence of a
bijection $\bN \leftrightarrow \DomWeights$ using his work on cells in
affine Weyl groups.  From the point of view of Harish-Chandra modules,
Vogan also conjectured a bijection between $\bN$ and $\DomWeights$.
Such a bijection has been established by Bezrukavnikov in two
preprints (the bijections in each preprint are conjecturally the same)
\cite{bezrukavnikov:tensor}, \cite{bezrukavnikov:quasi-exc}.
Bezrukavnikov's second bijection is closely related to Ostrik's
conjectural description of the bijection \cite{ostrik:equiv-k-theory}
(see also \cite{chm-ostrik:dist-inv}).  In the case of $G=GL(n, \C)$,
the first author \cite{achar:thesis} described an explicit
combinatorial bijection between $\bN$ and $\DomWeights$ from the
Harish-Chandra module perspective.  At present, it is not known how
all of these bijections are related (Bezrukavnikov's two candidates;
Ostrik's conjectural candidate; and the first author's candidate in
type $A$).
%although we make some comments on connections later in section 3.
In this paper, we work in the context of \cite{achar:thesis}, which we now 
review.

Let $K_G(\nilcone)$ denote the Grothendieck group of $G$-equivariant
coherent sheaves on $\nilcone$. On the one hand, $K_G(\nilcone)$ has a
natural basis indexed by elements of $\DomWeights$, denoted
$AJ(\lambda)$ in \cite{ostrik:equiv-k-theory} (but unlike in
\cite{ostrik:equiv-k-theory}, we are not utilizing the $\C^*$ action
on $\nilcone$ here).
%Essentially $AJ(\lambda)$ is a push forward of a line bundle 
%on the cotangent...
The algebraic global sections of $AJ(\lambda)$ are isomorphic as a
$G$-module to $\Ind_T^G \lambda$.  Thus, the space of global sections
of any element $\cF \in K_G(\nilcone)$ is given as a $G$-module by a
unique expression of the form
\begin{equation}\label{eqn:vogan-thm}
\Sections(\nilcone, \cF) = \sum_{\lambda \in \DomWeights}
m_\lambda \Ind_T^G \lambda,
\end{equation}
where the $m_\lambda \in \Z$, and $m_\lambda \ne 0$ for only finitely
many $\lambda$.  (This fact was communicated to us by Vogan.)

Let us now fix a pair $(e,\tau) \in \bN$.  We want to consider all
elements $\cF \in K_G(\nilcone)$ whose support is contained in
$\overline{\orb}_e$, and whose restriction to $\orb_e$ is the vector
bundle arising from $\tau$.  For each such $\cF$, there is at least
one $\lambda$ of maximal length occurring in the expression
(\ref{eqn:vogan-thm}) (where we have fixed a $W$-invariant
positive-definite symmetric bilinear form on the real span of
$\Weights$, so that we can speak of the length of a weight of $G$).
Define $\gamma: \bN \to
\DomWeights$ by
\[
\gamma(e,\tau) = \text{the \emph{smallest} such largest $\lambda$,
over all possible choices of $\cF$}.
\]
The following conjecture was made in \cite{achar:thesis}; moreover, it
was proved in the case of $G = GL(n,\C)$ in {\it op.~cit}.

\begin{conj}[\cite{achar:thesis}]
The map $\gamma$ is well-defined and is a bijection.  
Moreover, there is a basis
$\{M(e,\tau)\}$ for $K_G(\nilcone)$, indexed by $\bN$, such that
\begin{enumerate}
\item the support of $M(e,\tau)$ is $\overline{\orb}_e$,
\item the restriction of $M(e,\tau)$ to $\orb_e$ is the locally
free sheaf arising from $\tau$, and
\item there is an upper-triangular relationship between this basis and
the one indexed by $\DomWeights$:
\[
\Sections(\nilcone, M(e,\tau)) 
= \pm\Ind_T^G \gamma(e,\tau) + \sum_{\sz{\mu} < \sz{\tau}}
m_\mu \Ind_T^G \mu.
\]
\end{enumerate}
\end{conj}

In this paper, we study the problem of computing $\gamma$ when $\tau$
gives rise to a local system on $\orb_e$ (for semisimple groups, this
means that $\tau$ is trivial on the identity component of $G^e$).  Let
us denote by $\local$ the image of $\gamma$ when $\tau$ corresponds to
a local system.
%When $\tau$ is a local system, it seems likely that 
%$M(e, \tau)$ is just the direct image of the local system $\tau$
%on $\orb_e$ to the whole nilcone and so 
%$\Sections(\nilcone, M(e, \tau) = \Sections(\orb_e, \tau)$.
%This is consistent with Bezrukavnikov and Ostrik's work.
%(come later...)

Here is an outline of the paper.  In section 2, we compute $\local$
explicitly for $G=GL(n)$.  In section 3, we state a precise
conjecture, for general $G$, that the Dynkin weights for the Langlands
dual group of $G$ are a subset of $\local$.  In section 4, we show how
to associate to $\gamma(e, \tau)$ a sub-bundle of the cotangent bundle
of $G/B$.  Then for $GL(n)$, we are able to prove that when $e$ is in
a fixed nilpotent orbit, the cohomology of the sub-bundles (with
coefficients in the structure sheaf) are independent of the local
system $\tau$.

We thank David Vogan, Viktor Ostrik, and Roman Bezrukavnikov
for helpful conversations.
The second author was supported by NSF grant DMS-0070674.

\section{Computing $\gamma$ for a local system in $GL(n)$}

For a partition $\dd$ of $n$, let $\orb_\dd$ denote the nilpotent
orbit in $\gln(n)$ indexed by $\dd$.  If $\dd = [k_1^{a_1},
k_2^{a_2}, \ldots , k_l^{a_l}]$, then we know that
$\pi_1(\orb_\dd) \simeq \Z/c\Z$, where $c$ is the greatest common
divisor of the $k_i$'s \cite{col-mcg:nilp}.  Let $G^\dd$ denote the
isotropy group of an element in this orbit, and $\gred^\dd$ the
reductive part thereof.  Following the notation of
\cite{col-mcg:nilp}, we have
\[
\gred^\dd \simeq GL(a_1)_\Delta^{k_1} \times \cdots \times
GL(a_l)_\Delta^{k_l}.
\]
Here, $H_\Delta^m$ denotes the diagonal copy of $H$ in the product of
$m$ copies of itself.  Now, in the case of simply connected semisimple
groups, we identify local systems on an orbit with representations of the
component group of the centralizer of an element in the orbit.  
In the setting of
$GL(n)$, we will first produce a list of representations of
$\gred^\dd$, and then we will examine their restrictions to $SL(n)$ to
determine the correspondence with local systems.

Let $A_1 \in GL(a_1)$, $A_2 \in GL(a_2)$, etc.  The set of matrices
$(A_1, \ldots, A_l)$ determines an element $\tilde A$ of
\[
GL(a_1)_\Delta^{k_1} \times \cdots \times GL(a_l)_\Delta^{k_l};
\]
moreover, we have
\[
\det \tilde A = (\det A_1)^{k_1}\cdots(\det A_l)^{k_l}.
\]
Since every $k_i$ is divisible by $c$, we can write this as
\[
\det \tilde A = \big((\det A_1)^{k_1/c}\cdots(\det A_l)^{k_l/c}\big)^c.
\]
Let ``$\det^{1/c}$'' denote the character $\gred^\dd \to \C^\times$
given by
\[
(\det A_1)^{k_1/c}\cdots(\det A_l)^{k_l/c},
\]
and for any $p \in \Z$, let $\det^{p/c}$ denote the $p$-th power of
this map.

In $SL(n)$, the isotropy group $SL(n)^\dd$ is the
subgroup of $G^\dd$ consisting of matrices of determinant $1$.  The
function $\det^{1/c}$ takes only finitely many values when restricted
to this subgroup, namely, the $c$-th roots of unity.  Indeed, the
component group of $SL(n)^\dd$ is just $\Z/c\Z$, and the irreducible
representations of $\Z/c\Z$ are just powers of the function
$\det^{1/c}$.  We have, then, that the irreducible representations of
$\pi_1(\orb_\dd)$ come from 
\[
\det^{0/c}, \det^{1/c}, \ldots, \det^{(c-1)/c}.
\]

Now we need to compute $\gamma(\orb_\dd, \det^{p/c})$ for $0 \le p \le
c-1$.  Suppose that the dual partition of $\dd$ is $\dd\trans =
[j_1^{b_1}, \ldots, j_s^{b_s}]$.  Following the notation of
\cite{achar:thesis}, we let $L_\dd$ denote the Levi subgroup
\[
L_\dd = GL(j_1)^{b_1} \times \cdots \times GL(j_s)^{b_s},
\]
and we let $P_\dd$ be the standard (block-upper-triangular) parabolic
subgroup whose Levi factor is $L_\dd$.  Recall that every nilpotent
orbit in $GL(n)$ is Richardson; in particular, $\orb_\dd$ is Richardson
for $P_\dd$.  We fix a choice of $e \in \orb_\dd$
such that $\gred^{\dd} \subset L_\dd$ and $G^\dd \subset P_\dd$.
Explicitly, the embedding of $G^\dd$ in $P_\dd$ is such that the
factor $GL(a_i)^{k_i}_\Delta$ sits diagonally across the product of
the first $k_i$ factors of $L_\dd$.
Let
$\rho_\dd$ denote the half-sum of the positive roots of $L_\dd$.

We make the observation that any irreducible representation of $G^\dd$
must be trivial on its unipotent radical, hence is completely
determined by its restriction to $\gred^{\dd}$.  Therefore, it makes
sense to refer to irreducible representations of $G^\dd$ by their
highest weights.  We do this at various points below.

In \cite{achar:thesis}, the computational device of ``weight
diagrams'' was introduced and employed to carry out the computation of
$\gamma$.  For our present purposes, however, we are only concerned
with a small collection of representations for each orbit, and we do
not need to use the cumbersome weight diagrams.  Instead, we make use
of the following auxiliary result.

\begin{prop}[\cite{achar:thesis}, Claim~2.3.1]\label{prop:gamma}
Let $(\pi_\lambda, V_\lambda)$ denote the irreducible
$L_\dd$-representation of highest weight $\lambda$, regarded as a
$P_\dd$-representation by letting the unipotent part of $P_\dd$ act
trivially, and let $(\tau, W_\tau)$ be an irreducible representation
of $G^\dd$.  Suppose that $\tau$ occurs as a summand in
the restriction of $\pi_\lambda$ to $G^\dd$, and, moreover, that 
\[
\sz{\lambda + 2\rho_\dd}
\]
is minimal among all irreducible $P_\dd$-representations whose
restriction to $G^\dd$ contains $\tau$ as a summand.  Then
$\gamma(\orb_\dd, \tau) = \lambda + 2\rho_\dd$, made dominant for
$GL(n)$.
\end{prop}

%%should this go later
%To apply this proposition, we need to find a weight of $L_\dd$ whose
%restriction to $\gred^\dd$ is the (highest) weight of the character
%$\det^{p/c}$.  The weight of $\gred^\dd$ in question is
%\[
%\Bigg(\underbrace{\Big(\frac{pk_1}{c},\ldots,\frac{pk_1}{c}\Big)}
%_{GL(a_1)_\Delta^{k_1}},
%\ldots,\underbrace{\Big(\frac{pk_l}{c},\ldots,\frac{pk_l}{c}\Big)}
%_{GL(a_l)_\Delta^{k_l}}\Bigg).
%\]

%Let $\omega_p = \gamma(\orb_\dd, \det^{p/c})$.  
The next proposition describes $\gamma(\orb_\dd, \det^{p/c})$
explicitly in terms of the standard basis.
The remainder of the section will be devoted to establishing this
proposition.  In Section 4, we shall use this explicit
description to prove the main result. 

Let $\omega_p = \gamma(\orb_\dd, \det^{p/c})$.  Writing $\omega_0$ in
the standard basis, let ``block $\bB_a$'' refer to the collection of
coordinate positions which contain the entry $a$ in $\omega_0$, for
$a\in\Z$.  Let $\mu_k$ be the multiplicity of $k$ in $\dd\trans$.
(Thus, if there is some $j_i$ such that $j_i = k$, then $\mu_k = b_i$;
otherwise, $\mu_k = 0$.)  Note that every $b_i$ and every $\mu_k$ must
be a multiple of $c$.

\begin{prop}\label{prop:omega}
The length of block $\bB_a$ is $\sum_{k \geq 0}
\mu_{a+2k+1}$ when $a\geq0$ and 
the length of $\bB_{a}$ and $\bB_{-a}$ are the same.  
%c becomes d later
Moreover, if we write the length of block $\bB_a$ as $m_ac$, 
then $\omega_p$ is obtained
by replacing the first $m_ap$ entries in block $\bB_a$ with $a+1$.
\end{prop}

\begin{exam}
Consider the orbit labeled by $\dd = [6,3,3]$ in $GL(12)$.  Then
$c=3$, $\dd\trans = [3,3,3,1,1,1]$, and $\mu_3 =3$ and $\mu_1=3$.

We illustrate the preceding proposition by listing here all the
$\omega_p$.  The first sentence of the proposition describes
$\omega_0$ by giving the lengths of blocks, and the second sentence
tells us how to obtain the other $\omega_p$'s by modifying
$\omega_0$.
\begin{align*}
\omega_0 &= (\overbrace{2,2,2}^{\bB_2},
	\overbrace{0,0,0,0,0,0}^{\bB_0}, 
	\overbrace{-2,-2,-2}^{\bB_{-2}}) \\
\omega_1 &= (3,2,2,1,1,0,0,0,0,-1,-2,-2) \\
\omega_2 &= (3,3,2,1,1,1,1,0,0,-1,-1,-2)
\end{align*}
\end{exam}

\begin{proof}
Let $\lambda_p$ be the appropriate dominant weight of $L_\dd$ as
described in Proposition~\ref{prop:gamma}, such that $\omega_p$ is
just $\lambda_p + 2\rho_\dd$, made dominant for $GL(n)$.  In what
follows, we shall be careless about saying ``made dominant'' every
time; the reader should fill in those words wherever appropriate.

%Suppose that the dual partition to $\dd$ is $\dd\trans = [j_1^{b_1},
%\ldots, j_s^{b_s}]$, so that
%\[
%L_\dd = GL(j_1)^{b_1} \times \cdots \times GL(j_s)^{b_s},
%\]
%where each $b_i$ is a multiple of $c$.  
Let us begin with the
trivial representation of $\pi_1(\orb_\dd)$, namely $\det^{0/c}$.  Whatever
$\lambda_0$ is, it must be dominant for $L_\dd$, so that $\langle
\lambda_0, 2\rho_\dd \rangle \ge 0$.  Therefore
\[
\sz{\lambda_0 + 2\rho_\dd} = \sz{\lambda_0} + \sz{2\rho_\dd} +
2\langle \lambda_0, 2\rho_dd \rangle \ge \sz{2\rho_\dd}.
\]
Now $0$ is a weight of $L_\dd$ with the right restriction to $\gred^\dd$,
and taking $\lambda_0 = 0$ obviously minimizes $\sz{\lambda_0 +
2\rho_\dd}$ (the above inequality becomes an equality).  We therefore
have $\omega_0 = 2\rho_\dd$.  What does $2\rho_\dd$ look like?  For
each $GL(j_i)$ factor of $L_\dd$, we get a part of $2\rho_\dd$ that
looks like
\[
(j_i-1, j_i-3, \ldots, 1-j_i).
\]
Thus, in the total $2\rho_\dd$, a particular coordinate $a$ with $a \geq 0$ 
occurs
once for each factor $GL(j_i)$ with $j_i = a + 2k + 1$ for some $k$ with 
$k \geq 0$.
It follows that the length of block $\bB_a$ is precisely $\sum_{k \ge
0} \mu_{a+2k+1}$, as desired.  It is also clear that the length of 
block $\bB_a$ and block $\bB_{-a}$ are equal.
 
Next, we consider the case $p \ne 0$.  We will consider the first
factor $GL(a_1)_\Delta^{k_1}$ of $\gred^\dd$ individually; the other
factors would be treated identically.  The factor
$GL(a_1)_\Delta^{k_1}$ of $\gred^\dd$ sits diagonally across various
$GL(j_i)^{b_i}$ factors of $L_\dd$; indeed, it sits across $k_1$ of of
them.  Now, given a weight $\lambda$ of $L_\dd$, we obtain the
coordinates of the restriction $\lambda|_{\gred^\dd}$ by summing up
coordinates of $\lambda$ according to the diagonal embedding of the
factors of $\gred^\dd$ in $L_\dd$.  (See \cite{achar:thesis} for a
detailed account of how to restrict weights from $L_\dd$ to
$\gred^\dd$.)  In any case, we add up $k_1$ distinct coordinates of
$\lambda$ to obtain each coordinate of the $GL(a_1)_\Delta^{k_1}$ part
of the restriction of $\lambda$.  Now, the highest (and only) weight
of $\gred^\dd$ on the representation $\det^{p/c}$ is
\[
\Bigg(\underbrace{\Big(\frac{pk_1}{c},\ldots,\frac{pk_1}{c}\Big)}
_{GL(a_1)_\Delta^{k_1}},
\ldots,\underbrace{\Big(\frac{pk_l}{c},\ldots,\frac{pk_l}{c}\Big)}
_{GL(a_l)_\Delta^{k_l}}\Bigg);
\]
in particular, every coordinate in the $GL(a_1)_\Delta^{k_1}$ part of
the weight is equal to $pk_1/c$.  If we want $\sz{\lambda}$ to be
minimal, it is clear that the $k_1$ coordinates we add up to obtain
this coordinate should consist of $pk_1/c$ $1$'s and $(c-p)k_1/c$
$0$'s.

Repeating this argument for every factor of $\gred^\dd$, we see that
$\sz{\lambda}$ is minimized if, among all its coordinates, there are $pn/c$
$1$'s and $(c-p)n/c$ $0$'s.  But to compute $\gamma$, we need to
minimize $\sz{\lambda + 2\rho_\dd}$, not $\sz{\lambda}$.  We have
\[
\sz{\lambda + 2\rho_\dd} = \sz{\lambda} + \sz{2\rho_\dd} +
2\langle \lambda , 2\rho_\dd \rangle.
\]
So among possible $\lambda$'s of minimal size, we could try to choose
one so as to minimize $\langle \lambda, 2\rho_\dd \rangle$.  In fact,
we can arrange for the latter inner product to be $0$.

Among the factors $GL(j_i)^{b_i}$ of $L_\dd$, take the weight
$(1,\ldots,1)$ ({\it i.e.}, the determinant character) on $pb_i/c$ of
the factors, and $(0,\ldots,0)$ on the remaining $(c-p)b_i/c$ of
them.  Concatenating these weights together, we obtain a weight
$\lambda_p$ of $L_\dd$.  It is easy to see that $\lambda_p$ has the
right number of $0$'s and $1$'s to have the desired restriction to
$\gred^\dd$ as well as to be of minimal size.  Moreover, as promised, we
have that $\langle \lambda_p, 2\rho_\dd \rangle = 0$.

Therefore, $\omega_p = \lambda_p + 2\rho_\dd$.  This looks very
similar to $\omega_0$, except that in each block $\bB_a$, a proportion
$p/c$ of the coordinates have been increased by $1$.  Thus $\omega_p$
has precisely the desired form.
\end{proof}

%\begin{cor}\label{prop:trivrep}
%The weight $\omega_0 = \gamma(\orb_\dd, \det^{0/c})$
%is the Dynkin weight of the orbit corresponding to the partition $\dd\trans$.
%\end{cor}
%\begin{proof}
%We have showed that $\omega_0 = 2\rho_\dd$ (made dominant for $G$).
%But $2\rho_\dd$ is the Dynkin weight
%for the nilpotent orbit which is regular in $L_{\dd}$,
%which is exactly $\orb_{\dd\trans}$.
%\end{proof}

%We will verify the assertion of this proposition in reverse: it is
%claimed that if we take the weighted Dynkin diagram of
%$\orb_{\dd\trans}$, regard it as a weight, and rewrite it in the
%standard basis, we obtain precisely $\omega_0$.  Suppose that
%$\dd\trans = [j_1^{b_1}, \ldots, j_s^{b_s}]$.  To compute the
%corresponding weighted Dynkin diagram, we write out the numbers
%\begin{equation}\label{eqn:weight}
%j_1-1,j_1-3,\ldots,-j_1+1,\quad\cdot\,\cdot\,\cdot,\quad
%j_s-1,\ldots,-j_s+1,
%\end{equation}
%rewrite them in nonincreasing order, and compute all the differences
%of pairs of consecutive numbers in this list.  But if we stop a step
%early, before computing differences, what we have is the weight
%corresponding to the weighted Dynkin diagram of $\orb_{\dd\trans}$.
%Is this weight equal to $\omega_0$?  Let us go back one more step and
%consider the list of numbers in (\ref{eqn:weight}).  This list
%consists precisely of the entries of $2\rho_\dd$.  We saw in the
%course of the proof of Proposition~\ref{prop:omega} that this is
%precisely how $\omega_0$ is computed.

\section{A conjecture about Dynkin diagrams}

Let $\Langdual{G}$, $\Langdual{B}$, $\Langdual{T}$ be the data of the
Langlands dual group corresponding to $G$, $B$, $T$, respectively. Let
$\Langdual{\g}$, $\Langdual{\bo}$, $\Langdual{\h}$ denote the Lie
algebras of $\Langdual{G}$, $\Langdual{B}$, $\Langdual{T}$,
respectively. The weights of $T$ can be identified with the elements
$h \in \Langdual{\h}$ such that $\alpha\postcheck(h)$ is integral for
all coroots $\alpha\postcheck$ of $G$ (which are the roots of
$\Langdual{G}$).  This identification allows us to associate to a
nilpotent orbit $\Langdual{\orb}$ in $\Langdual{\g}$ a dominant weight
for $G$.  Namely, we can choose $e \in \Langdual{\orb}$ and let $e, h,
f$ span an $\Lie{sl}_2$-subalgebra of $\Langdual\g$ with $h \in
\Langdual{\h}$.  Then by Dynkin-Kostant theory, $h$ is well-defined up
to $W$-conjugacy, and by $\Lie{sl}_2$-theory, $h$ takes integral
values at the coroots of $G$.  Hence, $h$ determines an element of
$\DomWeights$.  We refer to the dominant weight of $G$ thus obtained
as the Dynkin weight of $\Langdual{\orb}$, and we denote by $\mathcal
D \subset \DomWeights$ the set of Dynkin weights of all nilpotent
orbits in $\Langdual{\g}$.
%See B-V for some motivation.

In the previous section, we saw that $\omega_0 = 2\rho_\dd$ (made
dominant for $GL(n)$).  This is nothing more than the Dynkin weight of
the nilpotent orbit $\orb_{\dd\trans}$ (we are identifying $G$ with
$\Langdual G$), since $\orb_{\dd\trans}$ intersects the Lie algebra of
$L_\dd$ in the regular orbit.  This result and a similar one for 
Richardson orbits in other
groups (along with calculations in groups of low rank) have led a
number of people to conjecture that $\mathcal D$ is a subset of
$\local$ (see \cite{chm-ostrik:dist-inv}).  We wish to state a precise
conjecture about how $\mathcal D$ sits inside of $\local$.  To this
end, we assign to each $\Langdual{\orb}$ in $\Langdual{\g}$ a finite
cover of $\orb= d(\Langdual{\orb})$, where $\orb$ is the nilpotent
orbit of $\g$ dual to $\Langdual{\orb}$ under Lusztig-Spaltenstein
duality.  Our conjecture essentially says that if we write our putative 
finite cover of $\orb$ as $G/ K$,
%(so $(G^e)^0 \subset H$)
then for some $\tau$ which is trivial on $K$,
$\gamma(\orb, \tau)$ equals the Dynkin weight of $\Langdual\orb$.

Let $A(\orb)$ denote the fundamental group of $\orb$ and let
$\bar{A}(\orb)$ denote Lusztig's canonical quotient of $A(\orb)$ (see
\cite{lusztig:reductive}, \cite{sommers:duality}).  Let
$\nilcone_{o,c}$ be the set of pairs $(\orb, C)$ consisting of a
nilpotent orbit $\orb \in \g$ and a conjugacy class $C \subset
A(\orb)$. We denote by $\dualnil_{o}$ the set of nilpotent orbits in
$\Langdual\g$.  In \cite{sommers:duality} a duality map $d:
\nilcone_{o,c} \to \dualnil_{o}$ is defined which extends
Lusztig-Spaltenstein duality. This map is surjective and the image of
an element $(\orb, C) \in \nilcone_{o,c}$, denoted $d_{(\orb, C)}$,
depends only on the image of the conjugacy class $C$ in
$\bar{A}(\orb)$.

Now given $\Langdual{\orb} \in \Langdual{\g}$, Proposition 13 of
\cite{sommers:duality} 
%and Remark 14
exhibits an explicit element $(\orb, C) \in \nilcone_{o,c}$
such that $d_{(\orb, C)}=\Langdual{\orb}$.  The orbit $\orb$ also
satisfies $\orb= d(\Langdual{\orb})$ (in particular, $\orb$ is special).  
Consider the image of $C$ in 
$\bar{A}(\orb)$, which we will also denote by $C$.  
We suspect that this conjugacy class coincides with one that Lusztig 
associates to $\Langdual{\orb}$
using the special piece for $\Langdual{\orb}$ (see Remark 14 in 
\cite{sommers:duality}).

The canonical quotient $\bar{A}(\orb)$ of $A(\orb)$ 
is always of the form 
$S_3, S_4, S_5$ or a product of copies of $S_2$.
Hence, it is possible to describe $\bar{A}(\orb)$ as a Coxeter group
of type $A$ and 
then to associate to each conjugacy class $C$ in $\bar{A}(\orb)$
a subgroup $H_C$ of $\bar{A}(\orb)$ which is well-defined up to 
conjugacy in $\bar{A}(\orb)$.
Lusztig did this for the exceptional groups in \cite{lusztig:notes}
and we now do it for the classical groups.

First we need to describe $\bar{A}(\orb)$ as a Coxeter group in the
classical groups (where $\bar{A}(\orb)$ is a product of copies of
$S_2$).  We use the description of $\bar{A}(\orb)$ in
\cite{sommers:duality}.  Let $\lambda = [\lambda_1^{a_1}, 
\lambda_2^{a_2}, \ldots, \lambda_k^{a_k}]$ 
be the partition corresponding to $\orb$ in the
appropriate classical group of type $B, C,$ or $D$.  Let $\mathcal{M}$
be the set of integers $m$ equal to some $\lambda_i$ such that
\begin{equation} 
\begin{aligned}
\lambda_i \text{ is odd and } \nu_i \text{ is odd} \text{ in type $B_n$} \\
\lambda_i \text{ is even and } \nu_i \text{ is even} \text{ in type $C_n$} \\
\lambda_i \text{ is odd and } \nu_i \text{ is even} \text{ in type $D_n$}  
\end{aligned}
\end{equation}
where $\nu_i = \sum^{i}_{j=1} a_j$. 
Then from section 5 of \cite{sommers:duality}, we know that the
elements of $\bar{A}(\orb)$ are indexed by subsets of $\mathcal{M}$ in
type $C$ and by subsets of $\mathcal{M}$ of even cardinality in types
$B$ and $D$.  In type $C$ we choose our set of simple reflections in
$\bar{A}(\orb)$ to correspond to subsets of $\mathcal{M}$ with a
single element.  In type $B$ and $D$ we choose our set of simple
reflections in $\bar{A}(\orb)$ to correspond to subsets $\{ a, b \}$
of $\mathcal{M}$ with $a > b$ and where no element of $\mathcal{M}$ is
both less than $a$ and greater than $b$.  Thus given a conjugacy class
$C$ of $\bar{A}(\orb)$ (which consists of a single element, $w$, since
the group is abelian), we can write $w$ minimally as a product of
simple reflections.  The simple reflections used are unique, and we
define $H_C$ to be the subgroup of $\bar{A}(\orb)$ generated by those
simple reflections.  Consider the surjection $\pi: G^{e} \to
\bar{A}(\orb)$ where $e \in \orb$ and define $K = \pi^{-1}(H_C)$ in
$G^{e}$.  We can now make our conjecture.

Given $\Langdual\orb$ in $\Langdual\g$, we have assigned a conjugacy
class $C$ in $\bar{A}(\orb)$ where $\orb = d(\Langdual\orb)$ and then
a subgroup $K$ in $G^e$ where $e \in \orb$.  Consider the finite cover
$\tilde{\orb}= G/ K$ of $\orb$.  Let $\C[\tilde{\orb}]$ denote the
global algebraic functions on $\tilde{\orb}$.  It is immediate that
$\C[\tilde{\orb}] = \sum \Gamma( \orb, \mathcal L)$ where the sum is
over the irreducible local systems $\mathcal L$ (counted with multiplicity)
which arise from the irreducible representations of
$A(\orb)$ appearing in $\Ind_{(G^e)^{0}}^{K} \C$, where $(G^e)^0$ is
the identity component of $G^e$. Hence (\ref{eqn:vogan-thm}) implies
that as a $G$-module $\C[\tilde{\orb}]$ can be written as 
$\sum_{\lambda \in \DomWeights} m_\lambda \Ind_T^G
\lambda$.  Let $\mu$ be a weight of largest length with $m_{\mu} \neq
0$.

\begin{conj}\label{conjecture:dynkin}
The weight $\mu$ is unique and is the Dynkin weight of
$\Langdual{\orb}$.
\end{conj}

We have verified the conjecture in a number of cases, although
a general proof is elusive at the moment.  

It seems likely that when $\tau$ gives rise to a local system, denoted
${\mathcal L}_{\tau}$, that $M(e, \tau)$ is just the direct image of
${\mathcal L}_{\tau}$ from $\orb_e$ to the whole nilpotent cone and so
$\Sections(\nilcone, M(e, \tau)) = \Sections(\orb_e, {\mathcal
L}_{\tau})$ (this would be consistent with Bezrukavnikov's and Ostrik's
work).  Since $\C[\tilde{\orb}] = \sum \Gamma( \orb, \mathcal L)$, 
Conjecture \ref{conjecture:dynkin} would
then state that the Dynkin weight of $\Langdual{\orb}$ occurs as
$\gamma(\orb, \tau)$ for some irreducible 
representation $\tau$ of $G^e$ which is trivial on $K$.

\begin{rmk} 
Specifying exactly what $\tau$ is 
seems to be more difficult.   For example, let 
$\Langdual\orb$ be the subregular orbit in type $B_n$.  Then
$\orb = d(\Langdual {\orb})$ 
is the smallest non-zero special orbit in type $C_n$.  This orbit
is Richardson, coming from the parabolic
subgroup with Levi factor of type $C_{n-1}$, and the parabolic 
subgroup gives rise to a $2$-fold cover of $\orb$.  Thus it is clear
that the Dynkin weight of $\Langdual{\orb}$ comes
from this $2$-fold cover of $\orb$ (which is, in fact, the one specified
by our conjecture) since $\Langdual{\orb}$ is regular in a
Levi factor of type $B_{n-1}$.
However, when $n$ is odd, the Dynkin weight will correspond
to the trivial local system on $\orb$, 
but when $n$ is even, the Dynkin weight will correspond to 
the non-trivial local system on $\orb$
(see the calculations in \cite{chm-ostrik:dist-inv}).
\end{rmk}

\section{Cohomology of the associated sub-bundles}

For an element $h \in \h$, we can define
a subspace $\nilrad_{h}$ of the nilradical $\nilrad$ of 
$\bo$ as follows.
We set
\[
\nilrad_{h} = \bigoplus_
{\substack{\text{$\alpha$ a positive root}\\
\alpha(h) \ge 2}} \Lie{g}_\alpha
\]
where $\Lie{g}_\alpha$ is the $\alpha$-eigenspace of the
root $\alpha$.
As in the previous section, 
we may identify $\DomWeights(\Langdual G)$ 
with a subset of $\h$.  Then for 
$\lambda \in \DomWeights(\Langdual G)$ we get 
a subspace $\nilrad_{\lambda}$ of $\nilrad$.
Our definition is motivated by the fact that
if $\lambda \in \h$ happens to be a Dynkin weight for a nilpotent
orbit $\orb \in \g$, then by work of McGovern,  
$\C[G \times_B \nilrad_{\lambda}] \simeq \C[\orb]$
\cite{mcgovern:regfn}
and moreover,
by work of Hinich and Panyushev,  
$H^i(G/B, S^j(\nilrad_{\lambda}^{*})) = 0$ for all $j \geq 0$ and $i>0$
\cite{hinich:van}, \cite{panyushev:van}.
Hence, it seems reasonable, especially given Conjecture
\ref{conjecture:dynkin}, to pick a general element 
$\lambda \in \local(\Langdual G)$ and study 
$H^i(G/B, S^j(\nilrad_{\lambda}^{*}))$.
Note that $\sum_{j \geq 0}
H^0(G/B, S^j(\nilrad_{\lambda}^{*})) \simeq \C[G \times_B \nilrad_{\lambda}]$.

We begin with the definition

\begin{defn}
Two $B$-representations $V, {\tilde V}$ 
are called $G/B$-equivalent if 
\[
H^i(G/B, S^j(V^*)) \simeq H^i(G/B, S^j({\tilde V}^*)).
\]
for all $i,j \geq 0$.
\end{defn}

Our main result is for $GL(n)$, and we identify $G$ with $\Langdual G$.
We hope in future work to say something interesting for other groups.

\begin{thm} \label{main}
Fix the partition $\dd$ and let $\nilrad_{\dd,p} = \nilrad_{\omega_p}$. 
For any $p$ the $B$-representations $\nilrad_{\dd,p}$ are all mutually 
$G/B$-equivalent to each other.
Thus, we have for any $p$ that
$$H^i(G/B, S^j(\nilrad_{\dd,p}^*)) \simeq
H^i(G/B, S^j(\nilrad_{\dd,0}^*))$$
and the latter equals $0$ if $i >0$ and
equals $\C^j[\orb_{\dd\trans}]$ if $i=0$
(the algebraic functions on the dual orbit of degree $2j$).
\end{thm}

The last sentence is just a re-formulation of the Hinich, McGovern, and
Panyushev results.   Before proving the theorem, we prove two lemmas which rely 
on the following proposition.

\begin{prop}\label{prop:demazure}
Let ${\tilde V} \subset V$ be representations of $B$ such that 
$V/{\tilde V} \simeq \C_{\mu}$
is a one-dimensional representation of $B$ corresponding to 
the character $\mu$.
Let $\alpha$ be a simple root, and let $P_\alpha$ be
the parabolic subgroup containing $B$ and having $\alpha$ as the only positive
root in its reductive part.  If $V$ extends
to a $P_\alpha$-representation and
$\langle\alpha\postcheck, \mu\rangle = -1$, then
$V, {\tilde V}$ 
are $G/B$-equivalent.
%\[
%H^i(G/B, S^j(V^*)) \simeq H^i(G/B, S^j({\tilde V}^*)).
%\]
%for all $i,j \geq 0$.
\end{prop}

\begin{proof}
Consider the exact sequence $0 \to {\tilde V} \to V \to \C_{\mu} \to 0$
and its dual $0 \to \C_{-\mu} \to V^* \to {\tilde V}^* \to 0$.
By Koszul, we have $0 \to S^{j-1}(V^*) \otimes \C_{-\mu}
\to S^j(V^*) \to S^j({\tilde V}^*) \to 0$ is exact for all $j \geq 0$.
By the lemma of Demazure \cite{demazure:bott-simple}, 
$H^*(G/B, S^{j-1}(V^*) \otimes \C_{-\mu}) = 0$ for all $j \geq 1$
since $V^*$ is $P_{\alpha}$-stable and
$\langle\alpha\postcheck, -\mu\rangle = 1$ (here our $B$ corresponds
to the positive roots, hence the difference with Demazure's convention).  
The result follows from the long exact sequence in cohomology.
\end{proof}

We need to introduce notation to 
describe the $B$-stable subspaces of $\nilrad$.
It is clear that if $\g_{\alpha}$ belongs to a 
$B$-stable subspace $U$ of $\nilrad$,
then so does $\g_{\beta}$ for all positive roots $\beta$
with $\alpha \preceq \beta$ (where $\preceq$
denotes the usual partial order on positive roots).
Hence it is enough to describe $U$ by the positive roots $\alpha$
such that $\g_{\alpha} \subset U$ and $\alpha$ is minimal among
all positive roots with this property.   
In this case, we say that $\alpha$ is minimal for $U$.

List the simple roots of $G$ as $\alpha_1, \dots, \alpha_{n-1}$.
Then any positive root of $G=GL(n)$ is of the form $\alpha_i +
\alpha_{i+1} + \cdots + \alpha_j$, which we denote by $[ i , j ]$.
We can express the usual partial order on the positive roots as $[ i ,
j ] \preceq [i', j']$ if and only if $i' \le i$ and $j \le j'$.  We
can then specify $U$ by its minimal positive roots, namely a
collection of intervals $[i,j ]$ such that for any two intervals $[
i,j ]$ and $[i' , j' ]$ with $i \leq i'$, we have $j \geq j'$. We will
say that $U$ is partially specified by the interval $[i,j]$ if the
root $[i,j]$ is minimal for $U$ (although there may be other minimal
roots).  Let us also say that $U$ is $i$-stable if $U$ is stable under
the action of the parabolic subgroup $P_{\alpha_i}$.

Let $U$ be a $B$-stable subspace of $\nilrad$ which is partially
specified by the interval $[a,b]$ and which is either
$(a\!-\!1)$-stable or $(b\!+\!1)$-stable.  Let $U'$ be the subspace of
$\nilrad$ which is specified by the same intervals as $U$ except that
$[a,b]$ is replaced by the two intervals $[a-1, b]$ and $[a, b+1]$.
Then $U$ and $U'$ are $G/B$-equivalent.  This is simply an application
of Proposition \ref{prop:demazure} where $\mu$ is the root $[a,b]$ and
$\alpha$ is either $\alpha_{a-1}$ or $\alpha_{b+1}$.  We refer to the
$G/B$-equivalence of $U$ and $U'$ as the \emph{basic move} for $a-1$
(respectively, for $b+1$).  We now state two lemmas which rely solely
on the basic move.

%\begin{lemma}
%Let $U$ be a subspace of $\nilrad$ 
%determined by  $[a,b]$ where $1 \leq a \leq b \leq n-1$.
%If $a \leq n-b$, 
%let $U'$ be the subspace of $\nilrad$ determined
%by $\{ [a-j+1, n-j] \ | \ 1 \leq j \leq a\}$ 
%If $a \geq n-b$,
%let $U'$ be the subspace of $\nilrad$ determined
%by $\{ [j, b+j-1] \ | \ 1 \leq j \leq n-b\}=
%\{ [n-b-j+1, n-j] \ | \ 1 \leq j \leq n-b\}$ 
%Then $R(U) = R(U')$. 
%\end{lemma}

\begin{lem}\label{lemma1}
Let $U$ be a subspace of $\nilrad$ which is partially specified
by an interval $[a,b]$ and such that $U$ is $i$-stable
for $b < i < d$ for some $d$.
Let $U'$ be the subspace of $\nilrad$ specified by the same intervals
as those defining $U$, but replacing $[a,b]$ with the collection 
of intervals
$$\{ [a-j+1, d-j] \ | \ 1 \leq j \leq d-b \}.$$ 
Then $U$ is $G/B$-equivalent to $U'$.
\end{lem}

\begin{proof}
Applying the basic move for $b+1$ replaces $[a,b]$ 
with the intervals $[a-1, b] \cup [a, b+1]$.
Now apply the basic move repeatedly on the right (for $i$
in the range $b+2 \leq i <d$)
to the interval $[a, b+1]$ and we are left with the two intervals
$[a-1, b] \cup [a, d-1]$.
The general result follows
by induction on $d-b$:  we apply the proposition
to the interval $[a-1, b]$ with $d$ replaced by $d-1$.
The base case $d-b=1$ is trivially true.
\end{proof}

\begin{lem}\label{lemma2}
Let $U$ be a subspace of $\nilrad$
partially specified by the intervals
$$[b_0, b_1] \cup [b_1+1,b_2] \cup \dots \cup 
[b_{l-1}+1, b_l] \cup [b_l+1, b_{l+1}]$$ 
and $$[b_1, b_2-1] \cup \dots \cup [b_{l-1}, b_l-1]$$
where $b_j \leq b_{j+1}-2$ for $1 \leq j \leq l-1$.
Assume that no interval of $U$ is of the form $[b_l, a_1]$ 
%or $[c, b_1-1]$; but
but that there are intervals of $U$ of the form $[b_0-1, a_2]$ 
and $[a_3, b_{l+1}+1]$.
Let $U'$ be the subspace of $\nilrad$ 
specified by the same intervals as $U$, except that 
$[b_0, b_1]$ is replaced by $[b_0, b_1-1]$
and $[b_l+1, b_{l+1}]$ is replaced by $[b_l, b_{l+1}]$.
Then $U$ is $G/B$-equivalent to $U'$.
A similar statement holds even if the we omit the interval $[b_l+1, b_{l+1}]$
or $[b_0, b_1]$ from the specification of $U$.
\end{lem}

\begin{proof}
Let $U_1$ be specified by the same intervals as $U$, except with the above 
intervals replaced by the intervals
$$[b_0, b_1] \cup [b_1+1,b_2-1] \cup [b_2+1, b_3-1]\cup \dots \cup [b_{l-1}+1, b_l-1] \cup [b_l+1, b_{l+1}]$$
Then $U_1$ is seen to be $G/B$-equivalent to $U$ by applying the
the basic move to the roots 
$[b_1+1,b_2-1], [b_2+1, b_3-1], \dots, [b_{l-1}+1, b_l-1]$ since
$U_1$ is stable for $b_2$, $b_3$, $\dots$ , $b_l$.  The stability for
$b_l$ follows from the assumption that no interval of $U$ is 
of the form $[b_l, a_1]$.

Let $U_2$ be specified by the same intervals as $U_1$,
except we replace the interval $[b_0, b_1]$ with $[b_0, b_1-1]$
and the interval $[b_l+1, b_{l+1}]$ with $[b_l, b_{l+1}]$.
Then $U_2$ is $G/B$-equivalent to $U_1$. 
This can be seen by applying the basic move
to $[b_0, b_1-1]$ (as $U_2$ is stable for $b_1$) and applying the
basic move (in reverse) to $[b_l, b_{l+1}]$ (as $U_1$ is $b_l$-stable).
Here we are using the fact that $U$ contains intervals
of the form $[b_0-1, a_2]$ and $[a_3, b_{l+1}+1]$.

The proof is completed by observing that 
$U_2$ is stable for $b_1$, $b_2$, $\dots$, $b_{l-1}$, 
so we can apply the basic move to the roots 
$[b_1+1,b_2-1], [b_2+1, b_3-1]$, $\dots$, $[b_{l-1}+1, b_l-1]$, arriving
at $U'$.  
\end{proof}

%\noindent {\bf Proof of Theorem.}
\begin{proof}[Proof of Theorem \ref{main}]
For an orbit corresponding to $\dd$ and a local system corresponding
to $p$, we have computed $\omega_p$ in Section 2.  We recall that $c$
is the g.c.d. of the parts of $\dd$ and $m_a$ was defined by the
equation $m_a c= \sum_{i \geq 0} \mu_{a+2i+1}$, where $\mu_i$ is the
multiplicity of $i$ as a part in $\dd\trans$.  Write $k$ for one less
than the largest part of $\dd\trans$, and set $s_i = \sum_{j \geq i}
m_j c$.

Then by Proposition \ref{prop:omega}, $\nilrad_{\dd,p} =
\nilrad_{\omega_p}$ is specified by the set of intervals $$I_i =
[s_{i+1} + p m_{i}, s_{i} + p m_{i-1}]$$ where $k \geq i \geq -k+1$.
%Let us set $a_i = s_{i+1} + p m_i$ to simplify notation,
%so $I_i = [a_i, a_{i-1}]$.
The difference between the left endpoint of $I_i$ and the right
endpoint of $I_{i+2}$ will be denoted $\Delta_i$.  So $$\Delta_i =
(c-p)m_{i+1} + p m_{i}.$$
%s_{i+1} + p m_{i} - (s_{i+2} + p m_{i+1})=
For $i \geq 1$ we have $\Delta_i \leq \Delta_{-i} \leq \Delta_{i-2}$
since $m_{-i}=m_i$ and $m_i \leq m_{i-2}$ when $i \geq 1$.

\medskip 

\emph{Step 1.}
Application of Lemma \ref{lemma1} to each of the intervals $I_i$ shows
that $\nilrad_{\dd,p}$ is $G/B$-equivalent to the subspace $W$ of
$\nilrad$ defined by the set of intervals $$I_{i,j} = [s_{i+1}+pm_i -
j+1, s_{i-1}+pm_{i-2} -j]$$ where $1 \leq j \leq \Delta_i$ if $k \geq
i \geq 1$ and where $1 \leq j \leq \Delta_{i-2}$ if $0 \geq i \geq
-k+1$.  In the former case, if $\Delta_{i} < j \leq \Delta_{i-2}$,
then the intervals $I_{i,j}$ are not minimal, so we may omit them from
our specification of $W$.
%Here, $a=s_{i+1}+pm_i$, $b=s_{i}+pm_{i-1}$, and 
%$d=s_{i-1}+pm_{i-2}$, so $d-b=\Delta_{i-2}$

\medskip

\emph{Step 2.}
Now for each $r$ in the range $k \geq r \geq 1$, starting with $r=k$
and working down to $r=1$, we will modify the intervals $I_{r+1,j}$
and $I_{-r+1, j}$ for $j > pm_{r+1}$ (in the cases where those
intervals are defined) to obtain a new subspace of $\nilrad$ which is
$G/B$-equivalent to $W$.

First, we will modify $I_{-r+1, j}$ for $\Delta_{r+1} < j \leq
\Delta_{-r-1}$.  Consider the intervals $$I_{r-1,j} \cup I_{r-3,j}
\cup \dots \cup I_{-r+5,j} \cup I_{-r+3, j}$$ for $\Delta_{-r-1} < j
\leq \Delta_{r-1}$
% (\leq \Delta_{-r+1})$
and the intervals $$I_{r-1,j} \cup I_{r-3,j} \cup \dots \cup
I_{-r+3,j} \cup I_{-r+1, j}$$ for $\Delta_{r+1} < j \leq
\Delta_{-r-1}$.  These yield a situation where we can apply Lemma
\ref{lemma2} a total number of $\Delta_{r-1} -\Delta_{-r-1} = p
(m_{r-1} - m_{-r-1})$ times to each of the intervals $I_{-r+1,j}$ for
$\Delta_{r+1} < j \leq \Delta_{-r-1}$.  This will replace
$I_{-r+1,j}=[s_{-r+2} + p m_{-r+1} - j+1, s_{-r} + p m_{-r-1}-j]$
%=[a_{-r+1} - j + 1, a_{-r-1}-j] 
with $[s_{-r+2} +pm_{-r-1} -j +1, s_{-r} +pm_{-r-1}-j]$.
%=[a_{-r+1}-j +1, a_{-r-1} + p(m_{r-1} - m_{-r-1})-j]$.
%Why can apply the second lemma:
%because for $\Delta_{r+1} < j \leq \Delta_{r-1}$ 
%there is no interval with a right endpoint equal to 
%$a_{r-1}-j$.

Second, we will modify both $I_{r+1,j}$ and $I_{-r+1, j}$ for
$pm_{r+1}<j \leq \Delta_{r+1}$ Consider the intervals $$I_{r-1,j} \cup
I_{r-3} \cup \dots \cup I_{-r+5,j} \cup I_{-r+3, j}$$ for
$\Delta_{r+1} < j \leq (\Delta_{r-1} -\Delta_{-r-1})+\Delta_{r+1}$ and
the intervals $$I_{r+1,j} \cup I_{r-1,j} \cup \dots \cup I_{-r+3,j}
\cup I_{-r+1,j}$$ for $pm_{r+1} < j \leq \Delta_{r+1}$.  We again
apply Lemma \ref{lemma2} a total number of $\Delta_{r-1}
-\Delta_{-r-1}$ times by modifying the pair of intervals $I_{r+1,j}
\cup I_{-r+1,j}$.  Then for $pm_{r+1} < j \leq \Delta_{r+1}$ we
replace $I_{r+1,j} \cup I_{-r+1,j}$ with $[s_{r+2} + pm_{r+1} -j
+1,s_{r} +pm_{r+1} -j] \cup [s_{-r+2} + pm_{-r-1} -j +1, s_{-r}
+pm_{-r-1}-j]$.
%=[a_{r+1} -j +1, a_{r-1} +p(m_{r+1} - m_{r-1})-j 
%=[a_{-r+1} + p(m_{-r-1} - m_{-r+1}) - j +1, a_{-r-1}-j]

We do this for all $r$ starting with $r=k$ and working backward to
$r=1$ (note we haven't done anything to the intervals $I_{1,j}$).

\medskip

\emph{Step 3.}
At this point our subspace of $\nilrad$ is specified by the following
intervals: $[s_{i} - l +1, s_{i-2} - l]$ for $i \geq 2$ with $1 \leq l
\leq (c-p)m_{i}$ and for $1 \geq i$ with $1 \leq l \leq (c-p)m_{i-2}$,
and those $I_{i,j}$ where $j \leq pm_i$ for $i \geq 1$ and $j \leq
pm_{i+2}$ for $i \leq 0$.
%Let us set $J_{i,l} = [s_{i+1} - l +1, s_{i-1} - l]$
%check again

Fix $r \geq 1$.  We may assume by induction on $r$ that our
subspace is $G/B$-equivalent to one partially specified by 
the following intervals (the case $r=1$ is already true):
$[s_{i} - l +1, s_{i-2} - l]$ for $1 \leq i \leq r$
and $1 \leq  l \leq cm_{i}$;
and $[s_{i} - l +1, s_{i-2} - l]$ 
for $0 \geq i \geq -r+2$ and $1 \leq l \leq cm_{i-2}$.
We want to show that the our subspace is $G/B$-equivalent
to one partially specified by the previous intervals 
together with the cases where $i=r+1$ and $i=-r+1$ in the 
above formulas.
% in the previous formulas.

Let $J_{i,l} = [s_{i+1} - l +1, s_{i-1} - l]$.
Consider the intervals
$$I_{r+1,j} \cup J_{r-2, j + (c-p)m_{r-1}} \cup \dots
\cup J_{-r+2 ,j + (c-p)m_{r-1}} \cup I_{-r+1,j}$$
for $1 \leq j \leq pm_{r+1}$ and the intervals
$$J_{r-2, l} \cup \dots \cup J_{-r+2,l}$$ 
for $1 \leq l \leq (c-p)m_{r-1}$.
%????  beware that $J_{r,l}$ kicks in at (c-p)m_{r+1}!!!!!!
We apply Lemma \ref{lemma2} a total number of 
%(in reverse??)
$(c-p)(m_{r-1} - m_{r+1})$ times to the
pair of intervals $I_{r+1,j} \cup I_{-r+1,j}$
for $1 \leq j \leq pm_{r+1}$.  We thus replace 
$I_{r+1,j} \cup I_{-r+1,j}$ with 
$[s_{r+2} + pm_{r+1} -j +1, s_{r-1}+(p-c)m_{r+1}-j] \cup
[s_{-r+1}+(p-c)m_{-r-1} -j+1, s_{-r}+pm_{-r-1}-j]$.
These intervals are just
$[s_{r+1} -l +1, s_{r-1}-l] \cup
[s_{-r+1} -l +1, s_{-r-1} -l]$
for $(c-p)m_{r+1} <l \leq cm_{r+1}$.
Hence, by induction on $r$ we see that our original
subspace is $G/B$-equivalent to 
the subspace specified by the intervals
$J_{i,l}$ where $1 \leq l \leq cm_{i+1}$ for $i \geq 1$
and where $1 \leq l \leq cm_{i-1}$ for $i \leq 0$. 
This subspace is independent of $p$ which completes the proof.
\end{proof}

\bibliography{pramod}
\bibliographystyle{pnaplain}
\end{document}